\documentclass[final,onefignum,onetabnum]{siuro250211}

\usepackage{lipsum}
\usepackage{amsfonts}
\usepackage{graphicx}
\usepackage{epstopdf}

 \usepackage{comment}

\usepackage{amsmath}
\usepackage{amssymb}
\usepackage{amscd}
\usepackage{bm}
\usepackage{url}
\usepackage{booktabs}
\usepackage[numbers]{natbib}
\usepackage{color,soul}
\usepackage{float}
\usepackage{amsmath, xparse}
\usepackage{commath}
\usepackage{multicol}
\usepackage{algorithm}
\usepackage{algpseudocode}
\usepackage{soul}
\usepackage{booktabs}
\usepackage{tabularx}
\ifpdf
  \DeclareGraphicsExtensions{.eps,.pdf,.png,.jpg}
\else
  \DeclareGraphicsExtensions{.eps}
\fi

\usepackage{enumitem}
\setlist[enumerate]{leftmargin=.5in}
\setlist[itemize]{leftmargin=.5in}

\newsiamremark{remark}{Remark}
\newsiamremark{hypothesis}{Hypothesis}
\crefname{hypothesis}{Hypothesis}{Hypotheses}
\newsiamthm{claim}{Claim}
\newsiamremark{fact}{Fact}
\crefname{fact}{Fact}{Facts}

\headers{Choice of Subspace for the Quasi-minimal Residual Method}{M. Hu}

\title{On the Choice of Subspace for the Quasi-minimal Residual Method for Linear Inverse Problems}

\author{Moshen Hu\thanks{Department of Mathematics, Emory University, Atlanta, GA (\email{moshen.hu@emory.edu}).}}

\dedication{\small\textit{Project advisor: Lucas Onisk \thanks{Department of Mathematics, Emory Univeristy, Atlanta, GA (\email{lonisk@emory.edu}).}}}

\begin{document}

\maketitle

\begin{abstract}
Inverse problems arise in various scientific and engineering applications, necessitating robust numerical methods for their solution. In this work, we consider the effectiveness of Krylov subspace iterative methods, including  GMRES,  QMR, and their range restricted variants for solving linear discrete ill-posed problems. We analyze the impact of subspace selection on solution quality. Our findings indicate that range restricted QMR can outperform standard QMR, and confirm the previously observed behavior that range restricted GMRES can be superior to conventional GMRES in terms of approximation efficacy. Notably, range restricted QMR demonstrates a key advantage over GMRES with respect to range restricted QMR's singular spectrum which can make the method less sensitive to errors that are naturally present making it particularly effective when the noise level in the problem is uncertain.
\end{abstract}

\section{Introduction} \label{sec:intro}
Inverse problems arise in various scientific and engineering applications where one seeks to determine unknown parameters from observed data. These problems are often ill-posed, which can mean that small perturbations in the available data can lead to significant deviations in the computed solutions. Inverse problems appear in fields such as medical imaging, geophysics, and signal processing, where recovering meaningful solutions from noisy or incomplete data is critical \cite{bertero1998inverseimaging}. To address the challenges of ill-posedness, regularization techniques are employed to stabilize the solution and mitigate the effects of noise error.

\subsection{Inverse problems}

Mathematically, inverse problems are often modeled as a linear system of equations 
\begin{equation} \label{ax=b}
    A x = b,
\end{equation}
where \( A \in \mathbb{R}^{n \times p} \) is a given matrix, \( x \in \mathbb{R}^{p} \) is the unknown solution, and \( b \in \mathbb{R}^{n} \) represents the observed data. In practice, the observed data  \(b\) are usually contaminated, often contains errors due to measurement noise, modeling inaccuracies, or other uncertainties. This can be expressed as

\[
b = b^{\text{exact}} + e
\]
where \( b^{\text{exact}} \) represents the unknown exact data, and \( e \) denotes the error term. Because approximating the solution of \eqref{ax=b} is often an ill-posed problem, directly solving the associated least-squares problem can propagate error into the solution, corrupting the recovery process. This necessitates the use of regularization or approximation techniques to obtain a stable and meaningful reconstruction \cite{Inverseproblem}.

\subsection{Linear discrete ill-posed problems}

An ill-posed problem lacks at least one of the following: existence, uniqueness, or stability of the solution. In this work, we focus on problems that lack stability, where small perturbations in data can lead to large variations in the solution. A well-posed problem is one that contains all three of the aforementioned properties whose definition is due to Hadamard\footnote{\url{https://en.wikipedia.org/wiki/Well-posed_problem\#cite_note-1}}.

When the discretized ill-posed problem at hand is a large, the coefficient matrix \( A \) is often ill-conditioned. This means that its condition number $\kappa(A)=\sigma_1 / \sigma_n$ 
is large (e.g, $10^3$), where \( \sigma_1 \) and \( \sigma_{\min\{n,p\}} \) are the largest and smallest singular values of \( A \), respectively. Herein, we consider the situation where the singular values of $A$ decay without significant gap so that small perturbations in \( b \) can cause large deviations in the least-squares solution \( x \), making direct solutions highly sensitive to noise. This instability makes traditional direct solution techniques, such as those utilizing the LU decomposition, unreliable for inverse problems \cite{IllposedBook}.

\subsection{Motivation and goal} In this work, we investigate the effect subspace selection has on the performance of the quasi-minimal residual (QMR) method for approximating the solution of linear discrete ill-posed problems. We then compare how this compares to the generalized minimal residual (GMRES) method under mis-estimated noise levels, with the goal of improving solution error in the presence of noise.

While GMRES and subspace restricted variants can be effective for solving non-symmetric large-scale linear systems, we found that its application to discrete ill-posed problems can be problematic when the noise level in the data is mis-estimated. In such cases, the lack of an accurate noise estimate may lead to over-iteration and amplification of noise, degrading the quality of the computed solution. The QMR method provides an alternative approach by modifying the subspace over which residual minimization is performed, offering better robustness when the noise level is underestimated. Works such as \cite{ellShiftGMRES_paper, SubspaceChoice, early_ellShift_paper, neuman2012algorithms} demonstrated that restricting the subspace (also referred to as \emph{range restriction}) where the solution is sought in GMRES can significantly improve the accuracy of the computed solution in discrete ill-posed problems. Motivated by these findings, we consider whether a similar range restriction strategy can enhance the performance of QMR, particularly in the context of non-symmetric ill-posed systems where noise levels are mis-estimated.

The remainder of this paper is organized as follows. Section 2 introduces Krylov subspaces and the GMRES and QMR methods. Section 3 introduces the range-restricted GMRES and range-restricted QMR methods where we provide the theoretical framework for the range-restricted iterative variants. Section 4 presents numerical experiments illustrating the behavior and performance of these methods. Section 5 offers concluding remarks.

\section{Background} \label{sec:background}

In this section, we introduce Krylov subspace methods, which play a fundamental role in the solution of large-scale linear systems and inverse problems. We begin by discussing the construction of Krylov subspaces and their significance in iterative methods. We then present the Arnoldi and Lanczos bi-orthogonalization procedures, which form the backbone of GMRES and QMR methods, respectively. The section ends with a brief discussion of the phenomenon of semiconvergence when considering iterative methods for inverse problems.

\subsection{Krylov subspaces}

Krylov subspace methods form a class of iterative techniques for approximating the solution of large linear systems of the form \eqref{ax=b}. Unlike direct solvers, such as LU decomposition, which require explicit matrix factorization and can be computationally expensive and memory-intensive for large-scale problems, Krylov subspace methods approach the associated least-squares problem of \eqref{ax=b} iteratively. These methods rely only on matrix-vector products with \( A \), and in some cases \( A^T \), making them well-suited for large, sparse, or structured linear operators. Moreover, due to their iterative nature, Krylov methods can exhibit inherent regularizing effects when applied to discrete ill-posed problems \cite{gazzola16}.

For a given matrix \( A \) and an initial residual vector \( r_0 = b - Ax_0 \), the Krylov subspace of dimension \( m \) is defined as:
\begin{equation} \label{original krylov space}
    K_m(A, r_0) = \text{span} \{r_0, Ar_0, A^2r_0, \dots, A^{m-1}r_0\}
\end{equation}
where $x_0$ is an initial guess.
The subspace \eqref{original krylov space} may be thought to encode increasing amounts of information about the solution, allowing the associated iterative method to update its approximations efficiently \cite{SubspaceChoice}. However, directly forming the columns of \eqref{original krylov space} produces a matrix that lacks orthogonality among its spanning vectors since applying \( A \) repeatedly to \( r_0 \) will converge toward the dominant eigenvector of \( A \), thus failing to capture an effective basis for the solution space. To overcome this issue, orthogonalization techniques such as Gram-Schmidt are employed within iterative techniques like the Arnoldi process to construct an orthonormal basis:
\begin{equation*}
    K_m(A, r_0) = \text{span} \{v_1, v_2, \dots, v_m\},
\end{equation*}
where each \( v_i \) is orthonormal to every other vector \( v_j \) for \( i \neq j \) in exact arithmetic \cite{SaadBook}. 


\subsection{Generalized minimal residual method} \label{sec:GMRES}
The GMRES method is an iterative solver designed to minimize the residual norm over the Krylov subspace \eqref{original krylov space} at each iteration. GMRES builds on the Arnoldi process to generate an orthonormal basis for the Krylov subspace and expresses the action of the matrix \( A \) through a smaller Hessenberg matrix. The following subsection presents the Arnoldi process, summarizes the key relations used in GMRES, and formulates the GMRES algorithm for solving large, non-symmetric square linear systems.

\subsubsection{Arnoldi process}

The Arnoldi process constructs an orthonormal basis for the Krylov subspace \( K_m(A, r_0) \)\eqref{original krylov space}, where \( A \in \mathbb{R}^{n \times n} \) is a square non-symmetric matrix and \( r_0 \) is the initial residual. Starting with \( v_1 = r_0 / \|r_0\| \), each new basis vector is obtained by applying \( A \) to the current vector and orthogonalizing the result against all previous basis vectors using the Gram-Schmidt procedure. Throughout, we will refer to \( \|\cdot\| \) as the 2-norm for vector norms.

At the $j^{th}$ iteration, the algorithm computes a matrix-vector product \( w_j = A v_j \), then projects \( w_j \) onto the existing basis \( \{v_1, \ldots, v_j\} \), and subtracts the projection to enforce orthogonality. The result is then normalized to produce \( v_{j+1} \). After $m$ iterations, the projection coefficients are stored in a matrix \( H_m \in \mathbb{R}^{(m+1) \times m} \), which takes an upper Hessenberg form due to the structure of the recurrence.

This iterative procedure yields the Arnoldi decomposition:
\begin{equation}\label{arnoldi}
A V_m = V_{m+1} H_m,
\end{equation}
where \( V_{m+1}\in \mathbb{R}^{n \times (m+1)} \) contains the orthonormal basis vectors, \( H_m \) is an upper Hessenberg matrix that encodes the action of \( A \) on the subspace, and \( V_{m}\in \mathbb{R}^{n \times m}\) are the first $m$ orthonormal columns of \( V_{m+1}\). The full process is given in Algorithm \eqref{alg:arnoldi}. For the solution of inverse problems, the number of Arnoldi steps is often small compared to the size of the linear system because the solution is often sufficiently approximated within a low-dimensional subspace \cite{lewis2009}.

\begin{algorithm}[H]
\caption{Arnoldi process}\label{alg:arnoldi}
\begin{algorithmic}[1]
\State Input: ${A} \in \mathbb{R}^{n \times n}$ and ${r_0} \in \mathbb{R}^n$
\State {Output:} ${V}_{m+1} \in \mathbb{R}^{n \times(m+1)}$ and ${H}_{ m} \in \mathbb{R}^{(m+1) \times m}$
\State Set ${v}_1 = {r_0} / \|{r_0}\|$
\For{$j=1, 2, \ldots, m$}
    \State Compute ${w}_j = {A} {v}_j$
    \For{$i=1, 2, \ldots, j$}
        \State $h_{i,j} = ({w}_j, {v}_i)$
        \State ${w}_j = {w}_j - h_{i,j} {v}_i$
    \EndFor
    \State $h_{j+1,j} = \|{w}_j\|$
    \If{$h_{j+1,j} = 0$}
        \State Stop
    \EndIf
    \State ${v}_{j+1} = {w}_j / h_{j+1,j}$
\EndFor
\end{algorithmic}
\end{algorithm}

\subsubsection{Summary of GMRES}

Using the Arnoldi relation \eqref{arnoldi} after $m$ steps, where \( V_m \) contains an orthonormal basis for the Krylov subspace \( K_m(A, r_0) \), GMRES computes the approximate solution of the least-squares problem associated with \eqref{ax=b} of the form
\begin{equation}\label{xes}
x_m = x_0 + V_m y_m.
\end{equation}
For the sake of exposition, we assume \( x_0 = 0 \). We note that the initial residual can be written as \( r_0 = b - A x_0 = \beta v_1 \), where \( v_1 \) is the first Arnoldi basis vector. Using \eqref{arnoldi}, the GMRES iterate \( x_m \) is given by \eqref{xes}. Here, the norm of the residual may be represented as: 
\[
\|b - A x_m\| = \|r_0 - A V_m y\| = \|\beta v_1 - V_{m+1} H_m y\|
\]
where we will denote $\|\cdot\|$ to represent the 2-norm throughout the rest of the work.
Here, \( \beta = \|r_0\| \) and \( e_1 \) represent the first principle axis vector in \( \mathbb{R}^{m+1} \). Since the columns of \(V_{m+1}\) are orthonormal and the Euclidean norm is invariant under orthogonal transformations, we may represent the small minimization problem at the $m^{th}$ step as follows 
\[
y_m = \arg \min_{y \in \mathbb{R}^m} \| \beta e_1 - H_m y \|
\] 
whose transformation from the original large-scale minimization problem into a much smaller least-squares problem of dimension \( (m+1) \times m \) can be solved efficiently. The complete algorithm can be found in detail in \cite{SaadBook} or \cite{GMRES_paper}.

\subsection{Quasi-minimal residual method} \label{sec:QMR}

The QMR method is an iterative solver designed for solving large, non-symmetric linear systems. It is built upon the Lanczos bi-orthogonalization process, which generates two bi-orthogonal bases for the Krylov subspaces associated with \( A \) and \( A^T \). These bases are then used to project the original problem onto a smaller subspace, where a quasi-minimal residual solution is computed. In contrast to methods like MINRES, which are tailored for symmetric systems and relies on the symmetric Lanczos relation, QMR is applicable to more general and non-symmetric problems.

\subsubsection{Lanczos Bi-orthogonalization for non-symmetric problems}

The Lanczos bi-orthog  $\allowbreak - \allowbreak$ onalization procedure extends the Lanczos process to handle non-symmetric matrices, ensuring pairwise orthogonality between two distinct sets of basis vectors \cite{SaadBook}. Unlike the Lanczos process, which relies on a single sequence of vectors to construct a tridiagonal matrix, the bi-orthogonalization approach generates two sequences of bi-orthogonal vectors—one for the matrix \( A \) and another for its transpose \( A^T \). These bi-orthogonal vectors form the basis for the reduction of \( A \) to a tridiagonal form.

The bi-orthogonal bases after $m$ steps for the two subspaces produced by the Lanczos bi-orthogonalization procedure are given by:
\begin{equation} \label{bi-subspacev}
\begin{split}
&{K}_m\left(A, v_1\right) = \operatorname{span}\left\{v_1, A v_1, \ldots, A^{m-1} v_1\right\} \\
&{K}_m\left(A^T, w_1\right) = \operatorname{span}\left\{w_1, A^T w_1, \ldots, \left(A^T\right)^{m-1} w_1\right\}.
\end{split}
\end{equation}
The initial vectors here given by \(v_1\) and \(w_1\) satisfy \((v_1, w_1) = 1\), ensuring proper scaling for the bi-orthogonalization process. Initialization is given by
\[
\beta_1 = \delta_1 = 0, \quad v_0 = w_0 = 0.
\]
For \(j = 1, 2, \ldots, m\), compute
\begin{equation} \label{lanczos_recurrence}
\alpha_j = (A v_j, w_j), \quad
\hat{v}_{j+1} = A v_j - \alpha_j v_j - \beta_j v_{j-1}, \quad
\hat{w}_{j+1} = A^T w_j - \alpha_j w_j - \delta_j w_{j-1},
\end{equation}
followed by
\[
\delta_{j+1} = \sqrt{|(\hat{v}_{j+1}, \hat{w}_{j+1})|}, \quad
\beta_{j+1} = \frac{(\hat{v}_{j+1}, \hat{w}_{j+1})}{\delta_{j+1}},
\]
and normalization:
\[
v_{j+1} = \frac{\hat{v}_{j+1}}{\delta_{j+1}}, \quad
w_{j+1} = \frac{\hat{w}_{j+1}}{\beta_{j+1}}.
\]
These steps build the tridiagonal matrix
\[
T_m =
\begin{bmatrix}
\alpha_1 & \beta_2 &        &        &        \\
\delta_2 & \alpha_2 & \beta_3 &        &        \\
        & \ddots   & \ddots   & \ddots &        \\
        &          & \delta_{m-1} & \alpha_{m-1} & \beta_m \\
        &          &        & \delta_m & \alpha_m
\end{bmatrix}
\]
with a three term recurrence relation given by the second expression of \eqref{lanczos_recurrence}.

Unlike the upper Hessenberg matrix from Arnoldi, the tridiagonal matrix \( T_m \) does not directly approximate the matrix \( A \), but it stores the coefficients of the linear combination of previously computed orthonormal basis vectors. Additionally, \(\{v_i\}_{i=1}^m\) is a basis for \(K_m(A, v_1)\), and \(\{w_i\}_{i=1}^m\) is a basis for \({K}_m(A^T, w_1)\). The following relations hold:
\begin{equation} \label{Lanczos_eqn1}
A V_m = V_m T_m + \delta_{m+1} v_{m+1} e_m^T \quad \& \quad A^T W_m = W_m T_m^T + \beta_{m+1} w_{m+1} e_m^T.
\end{equation}
In exact arithmetic, the product of $V_m$ and $W_m^T$ forms an identity matrix. The iterative process without a stopping criteria is given in Algorithm \eqref{alg:biortho} below.

\begin{algorithm}[H]
\caption{Lanczos bi-orthogonalization process}\label{alg:biortho}
\begin{algorithmic}[1]
\State Input: $A \in \mathbb{R}^{n \times n}$ and $b \in\mathbb{R} ^n$
\State Output: \({V}_{m+1} \in \mathbb{R}^{n \times (m+1)}, {T}_{m}\)
\State pick ($v_1,w_1$)=1
\For{\( j = 1, 2, \dots, m \)}
    \State \( \alpha_j = (A v_j, w_j) \)
    \State \( \hat{v}_{j+1} = A v_j - \alpha_j v_j - \beta_j v_{j-1} \), \( \hat{w}_{j+1} = A^T w_j - \alpha_j w_j - \delta_j w_{j-1} \)
    \State \( \delta_{j+1} = \sqrt{|(\hat{v}_{j+1}, \hat{w}_{j+1})|} \); If \( \delta_{j+1} = 0 \), {Stop}
    \State \( \beta_{j+1} = \frac{(\hat{v}_{j+1}, \hat{w}_{j+1})}{\delta_{j+1}} \)
    \State \( w_{j+1} = \hat{w}_{j+1} / \beta_{j+1} \), \( v_{j+1} = \hat{v}_{j+1} / \delta_{j+1} \)
    \State \( T_{j, j} = \alpha_j \), \( T_{j, j+1} = \beta_{j+1} \), \(  T_{j+1, j} = \delta_{j+1} \) (if \( j < m \))
    
\EndFor

\end{algorithmic}
\end{algorithm}

\subsubsection{Summary of QMR}

Having established the Lanczos bi-orthogonalization process, we now turn to its application in iterative solvers. One such method is the QMR algorithm, which uses the bi-orthogonal basis vectors generated to approximate the solution to non-symmetric linear least-squares problems. With the bi-orthogonal bases \eqref{bi-subspacev}, the QMR method utilizes the first of the two generated subspaces to approximate the solution of the least-squares problem associated with \eqref{ax=b}. It should be noted that QMR can be used to approximate the solution to a pair of coupled systems \cite{SaadBook}.

While GMRES maintains an orthonormal basis throughout its iterative process using the Arnoldi, the Lanczos bi-orthogonalization process generates two sets of  bi-orthogonal vectors. While these vectors are only pairwise orthogonal, we have observed that they provide a stable approximation framework for inverse problems. We can express an algebraically more efficient version of the first equation from \eqref{Lanczos_eqn1} by the following: 
\begin{equation}\label{bilanczos}
A V_m = V_{m+1} {T}_m,
\end{equation}
where \( V_m = [v_1, v_2, \dots, v_m] \) consists of one of the biorthogonal basis sets generated by the Lanczos biorthogonalization process, and \( V_{m+1} \) extends the subspace with an additional vector. The matrix \( {T}_m \) is an \( (m+1) \times m \) tridiagonal matrix which may be thought of as a projection matrix relating the two bi-orthogonal Krylov subspaces given by \eqref{bi-subspacev}.

In QMR, the \(m^{\text{th}}\) approximate solution \(x_m \in x_0 + K_m\) is given by
\[
x_m = x_0 + V_m y_m.
\]
The residual minimization may be written as
\[
\|b - A x_m\| = \min_{y \in \mathbb{R}^m} \|r_0 - A V_m y\| = \min_{y \in \mathbb{R}^m} \|V_{m+1} (\beta e_1 - T_m y)\|,
\]
where \( r_0 = \beta v_1 \) and \( \beta = \|r_0\| \). Unlike GMRES, where the basis vectors are orthonormal, in QMR, \(V_{m+1}\) consists of biorthogonal vectors, not orthonormal ones. Despite this, we assume that the norm is approximately invariant under the transformation $V_{m+1}$, allowing us to minimize the projected residual at $m^{th}$ iteration:
\[
 \min_{y \in\mathbb{R}^m} \|\beta e_1 - {T}_m y\|,
\]
where \(y \in \mathbb{R}^m\) is a vector of coefficients that needs to be determined. The vector \( \beta v_1 \) corresponds to \( \beta e_1 \), where \(e_1\) is the first principle unit vector in \(\mathbb{R}^{m+1}\), and the minimizer \(y_m\) can be computed efficiently by solving an \((m+1) \times m\) least-squares problem \cite{QMR_original}. The full process without stopping criterion is Algorithm \eqref{alg:qmr} below. 

\begin{algorithm}[H]
\caption{QMR algorithm}\label{alg:qmr}
\begin{algorithmic}[1]
\State {Input:}  \( A \in \mathbb{R}^{n \times n} \), \( b \in \mathbb{R}^n \)
\State {Output:} \( x_m \in \mathbb{R}^{n} \)
\State  \( r_0 = b - A x_0 \),  \( \beta = \|r_0\| \), \( v_1 = \frac{r_0}{\|r_0\| }  \)
\State Choose \( w_1 \) such that \( (w_1, v_1) = 1 \)
\State \( V = v_1 \), \( W = w_1 \),  \( {T} \)
\For{$j = 1, 2, \dots, m$}
    \State \( \alpha_j = (w_j, Av_j) \), \( \hat{r}_j = Av_j - \alpha_j v_j - \beta_j v_{j-1} \)
    \State \( \beta_{j+1} = \|\hat{r}_j\| \), if \( \beta_{j+1} < \epsilon \), break
    \State \( v_{j+1} = \hat{r}_j / \beta_{j+1} \)
    \State \( w_{j+1} = \hat{w}_{j+1} / \beta_{j+1} \), \( v_{j+1} = \hat{v}_{j+1} / \delta_{j+1} \)
    \State \( V = [V, v_{j+1}] \), \( W = [W, w_{j+1}] \), \( {T} \)
    \State  \( \min_y \| \beta e_1 - {T}_m y_m \| \) 
    \State \( x_m = x_0 + V_m y_m \)
\EndFor
\end{algorithmic}
\end{algorithm}

\subsection{Semiconvergence in Iterative Methods}\label{sec Semicon}

Iterative methods for solving ill-posed problems often exhibit a phenomenon known as semiconvergence \cite{Inverseproblem}. In the early iterations, the solution improves as the method captures the dominant components of the true solution, and the error decreases. However, as iterations continue, the influence of noise in the data becomes more pronounced, leading to a deterioration in solution quality. This behavior motivates regularization strategies such as the truncated singular value decomposition (TSVD), where only the components associated with large singular values are retained. Truncating the small, noise-amplifying singular values helps avoid the adverse effects of semiconvergence. 


\begin{figure}[H]
    \centering
    \includegraphics[width=0.50\textwidth]{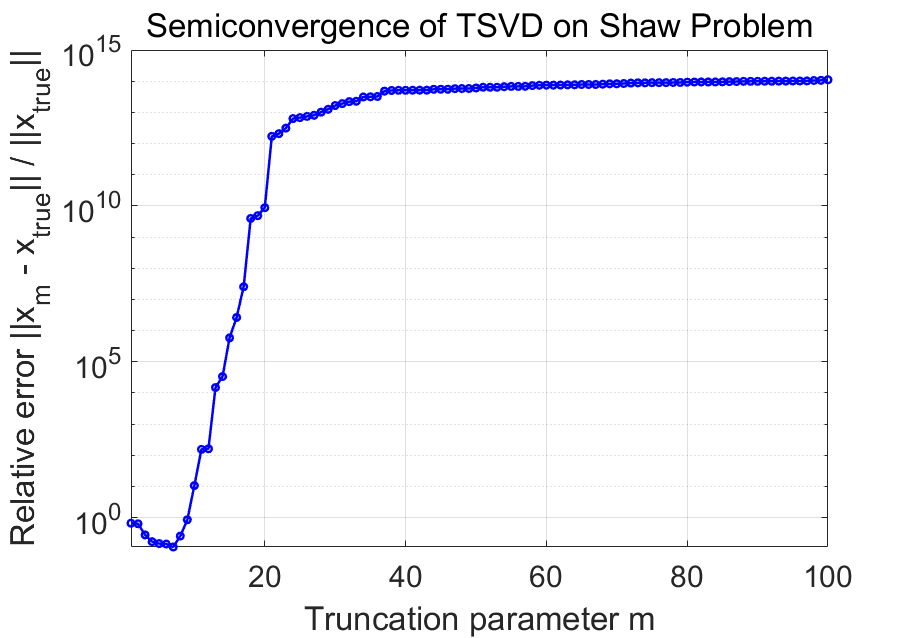}
    \caption{Semiconvergence behavior of the TSVD method applied to the Shaw\cite{Shaw_original} test problem with 1\% noise (see Section \ref{sec:results} for more details on experimental setup).}
    \label{fig:tsvd_semiconvergence}
\end{figure}

We summarize briefly how semiconvergence plays a role in discrete inverse problems. Suppose the noisy right-hand side is given by \( b = b^{\text{exact}} + e \). Then the exact solution satisfies
\[
x = A^{-1} b = A^{-1} b^{\text{exact}} + A^{-1} e.
\]
Assuming for the sake of argument that $A$ is square and non-singular, using the singular value decomposition (SVD) \( A = U \Sigma V^T \), we can express the solution as
\begin{equation} \label{nonsingular_soln}
x = \sum_{i=1}^{n} \frac{u_i^T b^{\text{exact}}}{\sigma_i} v_i + \sum_{i=1}^{n} \frac{u_i^T e}{\sigma_i} v_i.
\end{equation}
As the singular values \( \sigma_i \) decay rapidly for $A$, the reciprocals \( 1/\sigma_i \) become large, especially for large \( i \) which leads to amplification of even small noise components \( u_i^T e \). The second component of \eqref{nonsingular_soln} is often referred to as \emph{inverted noise}, which causes the solution to be dominated by noise if too many singular components of the solution are used \cite{Inverseproblem}.

The TSVD method demonstrates how selection of the dominant singular vectors can be used to stymie this propagation of error into the computed solution. Instead of including all \( n \) terms, TSVD retains only the first \( m \) terms corresponding to the $m$ largest singular values:
\[
x_m = \sum_{i=1}^{m} \frac{u_i^T b}{\sigma_i} v_i.
\]
Figure~\ref{fig:tsvd_semiconvergence} shows the relationship between the relative error and the truncation parameter \( m \). Initially, the reconstruction benefits from the inclusion of dominant components that capture meaningful features of the solution. However, as iteration parameter \( m \) increases and the smaller singular values fall below the noise level, additional terms begin to amplify noise rather than contribute useful information. This leads to a deterioration in solution quality, producing a characteristic U-shaped error curve. 

An analogous behavior appears in Krylov subspace methods such as GMRES and QMR, where the number of iterations acts as an implicit regularization parameter \cite{chung2024computational}. Stopping too soon results in under-regularization (i.e., over-smoothing), while excessive iterations introduce components aligned with small singular values, which are highly sensitive to noise. As shown in Figure \eqref{fig:tsvd_semiconvergence} of the numerical experiments, this results in deteriorating accuracy.

\section{Range restricted iterative methods} \label{sec:rangerestricted}

Krylov subspace methods such as GMRES and QMR work by approximating solutions within a subspace (or subspaces) generated by applying $A$, or, in the case of QMR, $A^T$, to the noisy right-hand side of an ill-conditioned system. While early iterations may provide reasonable approximations, prolonged iterations inevitably amplify errors into the computed solution, especially in ill-posed problems. To mitigate these issues a suitable termination criterion is applied. However, to generate a subspace (or subspaces) that is better aligned to certain problems inherent structures we can utilize range restricted Krylov methods. This section summarizes range restricted GMRES and derives the algorithm for QMR.

\subsection{Krylov subspace range restriction}
For range restricted Krylov methods, rather than generating iterates in the standard Krylov subspace
\begin{equation*}
K_m(A, b) = \text{span}\{b, A b, A^2 b, \dots, A^{m-1} b\},
\end{equation*}
we propose to restrict the iterates to the subspace
\begin{equation*}
K_m(A, A b) = \text{span}\{ A b, A^2 b, \dots, A^m b \}.
\end{equation*}
Here, we assume $x_0$ = 0, so that $r_0 = b$. When the signal or image we wish to recover is smooth, and the forward operator \( A \) acts as a smoothing operator, it can be more appropriate to search for a solution within a space whose vectors exhibit similar smoothness. By using the Krylov subspace generated from $Ab$ instead of $b$, we effectively begin in a space where the data has already been smoothed by the action of A.
Therefore, we can maintain focus on a more physically meaningful space that represents a better approximation to the true solution without amplifying unwanted noise \cite{early_ellShift_paper}. Furthermore, when \( A \) is a smoothing operator and the true underlying signal is smooth, the range of \( A \) tends to suppress noise present in \( b \), potentially damping the effect of noise on the degradation of the computed solution.



\subsection{The range restricted GMRES method} \label{sec:rrGMRES}

In range restricted GMRES, the shift parameter \( \ell \in \mathbb{N} \) defines how many times the space is shifted, effectively determining the depth of the transformation applied to the right-hand side. For instance, with \( \ell = 1 \), the method searches in the space spanned by \( A b \) and $A$, potentially filtering out noise in the null space of \( A \). With \( \ell = 2 \), the space further shifts to the subspace spanned by  \(A\) and \( A^2 b \), reinforcing a preference for smoother solutions. Unlike standard GMRES, which minimizes the residual norm in the Krylov subspace \(K_m(A, b) \), this modification ensures iterates remain within a shifted subspace, leveraging the properties of \( A \) to mitigate noise amplification \cite{early_ellShift_paper}.

Similar to standard GMRES, range-restricted GMRES also requires an orthonormal basis. We recall the Arnoldi decomposition \eqref{arnoldi} where $V_{m+1}$ has orthonormal columns and $H_{m}$ is an upper Hessenberg matrix. Using the QR factorization of $H_{m}$,
\begin{equation}\label{H=QR}
    H_{m} = Q^{(1)}_{m+1} R^{(1)}_m,
\end{equation}
where $Q^{(1)}_{m+1}\in\mathbb{R}^{m\times m}$ is an orthogonal matrix and $R^{(1)}_m\in\mathbb{R}^{(m+1) \times m}$ is upper triangular, we define:
\begin{equation*}
    W^{(1)}_m = V_{m+1} Q^{(1)}_{m+1}.
\end{equation*}
From \eqref{H=QR}and \eqref{arnoldi}, it follows that 

\begin{equation}\label{W=AVR}
    W^{(1)}_m = V_{m+1} Q^{(1)}_{m+1} = A V_m (R^{(1)}_m)^{-1}.
\end{equation}
Using the Arnoldi relation \eqref{arnoldi} gives:
\begin{equation*}
    A V_m (R^{(1)}_m)^{-1} = V_{m+1} Q^{(1)}_{m+1},
\end{equation*}
implying that the columns of $W^{(1)}_m$ span $K_m(A, A b)$. Since $V_m$ already spans $K_m(A, b)$, we obtain $A V_m = AK_m(A, b) =K_m(A, A b)$.

In the case of 2 shifts, we present a brief overview of the Arnoldi correspondent to the 2-shifted GMRES method. The matrix \( {W}_m^{(2)} \) is defined as the first \( m \) columns of \( {V}_{m+2} {Q}_{m+2}^{(2)} \). In the case of a 2-shift, the relation  

\begin{equation}\label{W2}
{W}_m^{(2)} = {A} {W}_m^{(1)}\left({R}_m^{(2)}\right)^{-1}
\end{equation}
follows from \eqref{W=AVR} and ensures that the column space of \( {W}_m^{(2)} \) corresponds to the shifted Krylov subspace \(K_m({A}, {A}^2 {b}) \), effectively incorporating the second shift into the iterative framework.
To generalize this to $\ell$-shifts, we recursively define $W_m^{(\ell)} = A W_m^{(\ell-1)} (R_m^{(\ell)})^{-1}$, which ensures that the columns of $W_m^{(\ell)}$ span $K_m(A, A^\ell b)$. Full details may be found in \cite{ellShiftGMRES_paper}.

In the case of more than one shift, successive QR factorizations play a crucial role in the algorithm. They ensure that the columns remain orthonormal, forming a well-conditioned basis for the subspace, while also guaranteeing that the span of the vectors used in the computation accurately represents the restricted Krylov subspace of interest. Additionally, QR factorizations provide a natural algorithmic framework for implementing these numerical methods effectively \cite{neuman2012algorithms}. In the case of multiple shifts, each step involves computing a new factorization:
\begin{equation} \label{succqr}
    H_{m+\ell+1, m+\ell} Q_{m+\ell, m}^{(\ell)} = Q_{m+\ell+1}^{(\ell+1)} R_{m+\ell+1, m}^{(\ell+1)},
\end{equation}
where $Q_{m+\ell+1}^{(\ell+1)}\in\mathbb{R}^{(m+\ell+1)\times (m+\ell+1)}$ is an orthogonal matrix and $R_{m+\ell+1, m}^{(\ell+1)}\in\mathbb{R}^{(m+l+1)\times m}$ is upper triangular. This stepwise factorization ensures that each transformation aligns the new basis with the shifted Krylov subspace and prevents loss of orthogonality due to rounding errors. The recursive structure of these QR factorizations allows efficient computation while preserving the structure of the projected system.

The minimization problem is then formulated as:

\begin{equation*}
    \min_{x\in K_m(A, A^\ell b)} \| A x - b \| = \min_{y \in\mathbb{R}^m}\|AW_m^{(\ell)}y-b\|.
\end{equation*}
Expressing the solution in terms of $y$,
\begin{equation*}
    \min_y \| A V_{m+\ell} Q^{(\ell)}_{m+\ell} y - b \|,
\end{equation*}
and using equation \eqref{W2} and further QR factorizations \eqref{succqr}, we may rewrite the previous expression as follows:
\begin{equation*}
    \begin{split}
        \min_y & \| V_{m+\ell+1} Q^{(\ell+1)}_{m+\ell+1} R^{(\ell+1)}_m y - b \|.
    \end{split}
\end{equation*}

Since the 2-norm is preserved under orthogonal transformations and the first column of $V_{m+\ell+1}$ is $b/\|b\|$, the problem reduces:
\begin{equation*}
    \min_y \| R^{(\ell+1)}_m y - \beta (Q^{(\ell+1)}_m)^T e_1 \|.
\end{equation*}
This reduced system is relative easy to solve, since \( R^{(\ell+1)}_m \) is upper triangular, allowing for efficient back-substitution. The final solution is given by  

\[
x_m^{(\ell)} = W_m^{(\ell)} y_m.
\]

The full process without  termination criterion is Algorithm \eqref{alg:rrgmres} below. The range restricted GMRES method provides an approach to solving ill-conditioned linear systems while ensuring iterates remain in $K_m(A, A^\ell b)$. This adjustment retains the efficiency of GMRES, making it particularly useful for solving ill-posed problems
\cite{ellShiftGMRES_paper}.

\begin{algorithm}[]
\caption{Range restricted GMRES ($\ell \geq 1$) }\label{alg:rrgmres}
\begin{algorithmic}[1]
\State {Input:} ${A} \in \mathbb{R}^{n \times n}, {b} \in \mathbb{R}^n$, and $\ell \in\{1,2,3, \ldots\}$
\State {Output}${x}_m^{(\ell)} \in \mathbb{R}^n$
\State ${v}_1 = {b} / \|{b}\|$ and ${x}_0^{(\ell)} = 0$
\For{$i = 1,2, \ldots, \ell$}
    \State Compute $\ell$ steps of Arnoldi $AV_{i} = V_{i+1} H_{i+1, i}$
\EndFor
\For{$m = 1,2, \ldots$}
    \State Compute next Arnoldi step ${A} {V}_{\ell+m} = {V}_{\ell+m+1} {H}_{\ell+m+1, \ell+m}$
    \State Compute QR factorization
    $\left[{Q}_{m+1}^{(1)},{R}_{m+1,m}^{(1)}\right] = {H}_{m+1, m}$
    \For{$j = 1,2, \ldots, \ell$}
        \State 
        Compute QR factorization$\left[{Q}_{j+m+1}^{(j+1)}, {R}_{j+m+1,m}^{(j+1)}\right] = {H}_{j+m+1, j+m} {Q}_{j+m, m}^{(j)}$
        
    \EndFor
    \State 
    $\min_y \left\|{R}_{\ell+m+1, m}^{(\ell+1)}{y} - \|{b}^\delta\| \left({Q}_{\ell+m+1}^{(\ell+1)}\right)^T{e}_1\right\|$
    
    \State ${x}_m^{(\ell)} = {V}_{\ell+m} {Q}_{\ell+m, m}^{(\ell)} {y}_m^{(\ell)}$
\EndFor
\end{algorithmic}
\end{algorithm}

\subsection{The Range Restricted QMR Method} \label{sec:rrQMR}

The range restricted QMR method extends the idea of Krylov subspace shifting to the QMR framework. Instead of searching for the solution using \( K_m(A, b) \) and \( K_m(A^T, b) \), the method utilizes the shifted spaces \( K_m(A, A^\ell b) \) and \( K_m(A^T,\left(A^T\right)^\ell b) \), in a similar vein to range restricted GMRES. While QMR relies on the Lanczos bi-orthogonalization process and requires access to both \( A \) and \( A^T \), the subspace restriction similarly aims to suppress noise by confining the iterates to a smoother, more stable component of the range of \( A \) while approximating the solution of \eqref{ax=b}.

Starting from the bi-orthogonal Lanczos decomposition \eqref{bilanczos}, one-time QR factorization is applied $T_{m} = Q^{(1)}_{m+1} R^{(1)}_m$. Then, a sequence of QR factorizations is recursively applied to \( T_{m+\ell}Q^{(j)}_{m+\ell} \) from $j=1,\dots,\ell$, producing a nested set of orthonormal bases that define the range-restricted subspace. Similarly to range restricted GMRES, we define
\begin{equation} \label{W=AWR}
W_m^{(\ell)} = A V_{m+\ell} (R_{m+\ell}^{(\ell)})^{-1}
\end{equation}
where \( W_m^{(\ell)} \) spans the subspace associated with \( K_m(A, A^\ell b) \) and \( K_m\left(A^T, \left(A^T\right)^\ell b\right) \). 

We present a brief overview of the 2-shifted QMR method. Similar to range restricted GMRES, the matrix \( {W}_m^{(2)} \) is defined as the first \( m \) columns of \( {V}_{m+2} {Q}_{m+2}^{(2)} \). In the case of a 2-shift, the relation  

\begin{equation*}
{W}_m^{(2)} = {A} {W}_m^{(1)}\left({R}_m^{(2)}\right)^{-1}
\end{equation*}
follows from \eqref{W=AWR} and ensures that the column space of \( {W}_m^{(2)} \) corresponds to the shifted joint Krylov subspace. The range restricted QMR iterate is then sought in the form
\[
x_m = x_0 + W_m^{(\ell)} y_m,
\]
with \( y_m \in \mathbb{R}^m \) minimizing the residual:
\[
y_m = \arg \min_y \| A W_m^{(\ell)} y - b \|.
\] As in the GMRES case, orthogonal transformations reduce this problem to solving an upper-triangular least-squares system:
\[
\min_y \| R^{(\ell+1)}_m y - \beta (Q^{(\ell+1)}_m)^T e_1 \|.
\]
The full algorithm without termination criterion is outlined in Algorithm 3.2 below.

\begin{algorithm}[H]
\caption{Range restricted QMR ($\ell \geq 1$)}
\begin{algorithmic}[1]
\State Input: ${A} \in \mathbb{R}^{n \times n}, {b}\in \mathbb{R}^n$, and $\ell \in\{1,2,3, \ldots\}$
\State Output: ${x}_m^{(\ell)} \in \mathbb{R}^n$
\State ${v}_1 = {b} / \|{b}\|$ and ${x}_0^{(\ell)} = 0$
\For{$i = 1,2, \ldots, \ell$}
    \State Compute $\ell$ steps of Arnoldi $A V_{i} = V_{i+1} T_{i+1, i}$
\EndFor
\For{$m = 1,2, \ldots$}
    \State Compute next Arnoldi step ${A} {V}_{\ell+m} = {V}_{\ell+m+1} {T}_{\ell+m+1, \ell+m}$
    \State  
    Compute QR factorization $\left[{Q}_{m+1}^{(1)},{R}_{m+1,m}^{(1)}\right] = {T}_{m+1, m}$
    \For{$j = 1,2, \ldots,\ell$}
        \State Compute QR factorization$\left[{Q}_{j+m+1}^{(j+1)},{R}_{j+m+1,m}^{(j+1)}\right] = {T}_{j+m+1, j+m} {Q}_{j+m, m}^{(j)}$
    \EndFor
    \State$\min_y \left\|{R}_{\ell+m+1, m}^{(l+1)}{y} - \|{b}\| \left({Q}_{\ell+m+1}^{(\ell+1)}\right)^T{e}_1\right\|$
    \State ${x}_m^{(\ell)} = {V}_{\ell+m} {Q}_{\ell+m, m}^{(\ell)} {y}_m^{(\ell)}$
\EndFor
\end{algorithmic}
\end{algorithm}

\section{Numerical results} \label{sec:results}

This section presents numerical experiments that demonstrate the effectiveness of range restricted QMR in solving ill-posed linear systems. We compare the performance of QMR and GMRES, as well as range restricted GMRES and QMR. The experiments cover a range of test problems, including one-dimensional and two-dimensional inverse problems to illustrate how restricting the solution space can enhance numerical robustness. By analyzing relative error, residual norms, and singular value behavior, we provide insights into the potential benefits of range restricted QMR in applications.
\subsubsection{Preliminaries}
In this section, we outline the key evaluation criteria and test problems used in our numerical experiments. We describe the discrepancy principle, which will serve as the algorithmic stopping criterion, metrics for solution evaluation, and the specific test cases considered in both one-dimensional and two-dimensional settings. To simulate noise, we generate a vector \( e \in \mathbb{R}^{1000} \) with entries drawn from a normal distribution with zero mean. The perturbed right-hand side is then given by \( b = b^{exact} + e \), where the noise vector \( e \) is scaled to achieve a prescribed noise level defined by
\[
v = 100 \cdot \frac{\|e\|}{\|b\|}.
\]
We will refer to \emph{v} as the noise level. We consider noise levels of 5\%, 1\%, 0.5\%, and 0.1\%. These settings establish the framework for the numerical results presented in the following sections.

\subsubsection{Termination criterion: discrepancy principle}

To avoid overfitting noisy data in ill-posed problems, iterative methods should not minimize the residual indefinitely. The discrepancy principle provides a practical stopping criterion: terminate the iteration when the residual norm falls below a threshold proportional to the noise level. Specifically, if the noise \( e \) satisfies \( \|e\| \leq \epsilon \), then the iteration is stopped once
\[
\|A x_m - b\| \leq \eta \epsilon,
\]
where \( \eta > 1 \) is a user-defined safety factor, typically chosen close to one (e.g., \( \eta = 1.01 \)) to account for noise estimation uncertainty.

This criterion ensures that the approximate solution \( x_m \) fits the data up to the level of noise, thereby preventing semiconvergence-related deterioration \cite{deblurringBook}.

\subsubsection{Solution evaluation}

To assess the performance of reconstruction algorithms, we use two standard metrics: relative residual and relative error. 

The relative residual evaluates how well the computed solution \( x \) satisfies the linear system and is defined by
\[
\frac{\|b - A x\|}{\|b^{exact}\|}.
\]
This metric reflects the fidelity of the solution to the observed data. However, in ill-posed problems, a small residual may not imply a good solution due to the presence of noise, as discussed in the semi-convergence section. \eqref{sec Semicon}

The relative error measures the accuracy of the computed solution relative to the true solution \( x_{true} \):
\[
\frac{\|x_{true} - x\|}{\|x_{true}\|}.
\]
This provides a direct assessment of reconstruction quality, though it is typically computable only for synthetic or test problems where \( x_{true} \) is known.

\subsubsection{Test problems: 1D and 2D cases}

The numerical experiments presented in this paper involve both one-dimensional (1D) and two-dimensional (2D) inverse problems arising from integral equations and image deblurring.

The 1D problems considered are based on the discretization of Fredholm integral equations of the first kind, using numerical methods such as the Nystr\"om method with the trapezoidal rule. These problems include the Phillips and Shaw problems. The Phillips problem is formulated from a classical integral equation introduced by D.~L.~Phillips \cite{Phillips_original}, with a kernel, solution, and right-hand side defined through a function incorporating cosine terms and piecewise structure. The Shaw problem models 1D image restoration and is defined by a kernel involving sine and cosine terms. The exact solution consists of a sum of Gaussian functions. This problem simulates typical blurring behavior encountered in signal processing applications.

The 2D problems involve image deblurring, in which a blurred image is modeled as the result of applying a linear blurring operator to an original image. The goal is to reconstruct the original image from noisy, blurred observations. These problems are formulated as large-scale linear inverse problems, where the blurring operator is discretized as a matrix \( A \), often defined by a Gaussian point spread function (PSF). Additional blur types such as motion blur, rotational blur, and atmospheric turbulence may also be modeled. We use IR Tools \cite{IRtools} to generate and test these problems.

\subsection{Numerical experiments}

In this section, we apply the preliminaries and problem formulations discussed earlier to evaluate the performance of the algorithms investigated in the previous Sections. Through numerical experiments, we analyze their stability, accuracy, and convergence behavior.

\subsection{Shifted VS. Non-Shifted: GMRES and QMR}

\begin{figure}[H]
    \centering
    \includegraphics[width=0.32\linewidth]{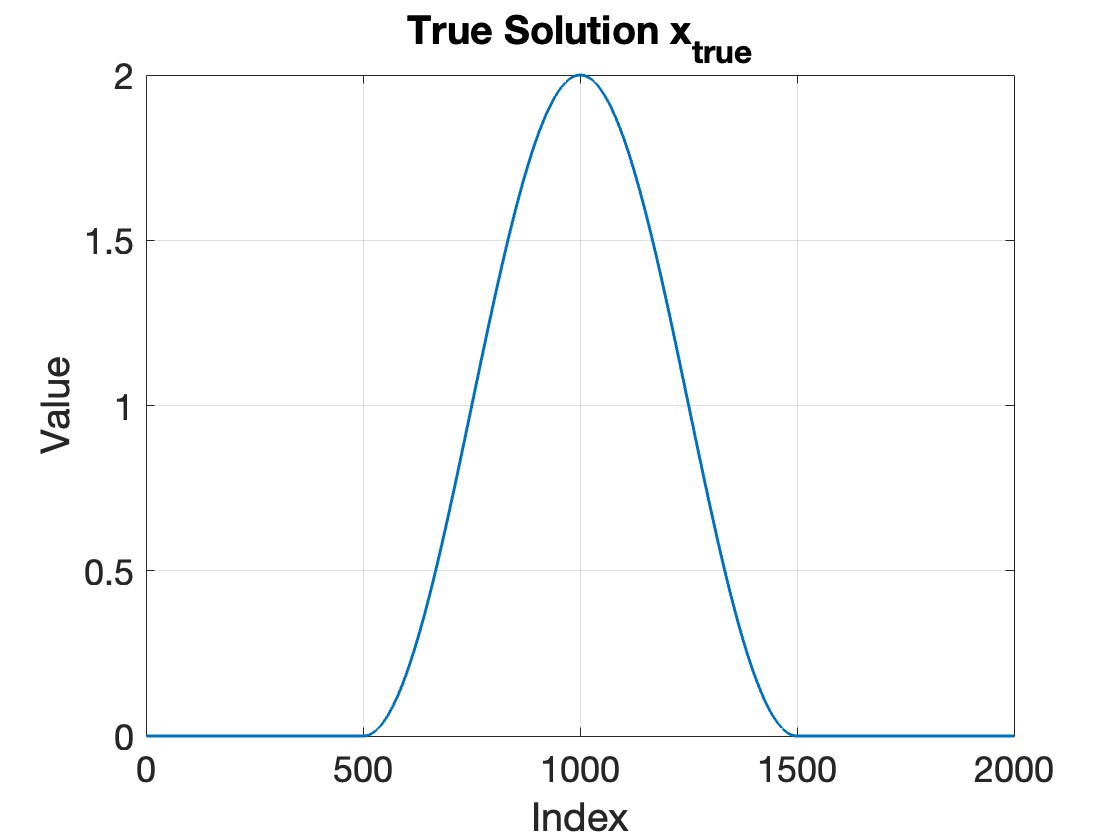}
    \includegraphics[width=0.32\linewidth]{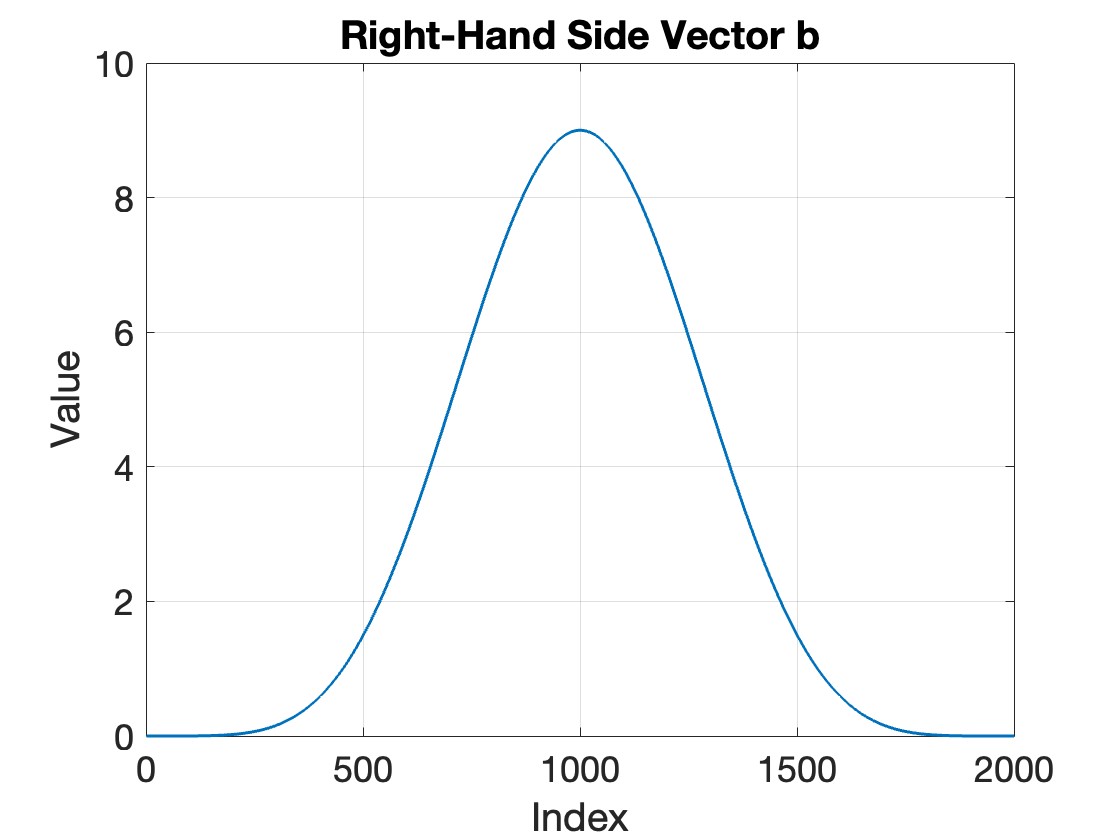}
    \includegraphics[width=0.32\linewidth]{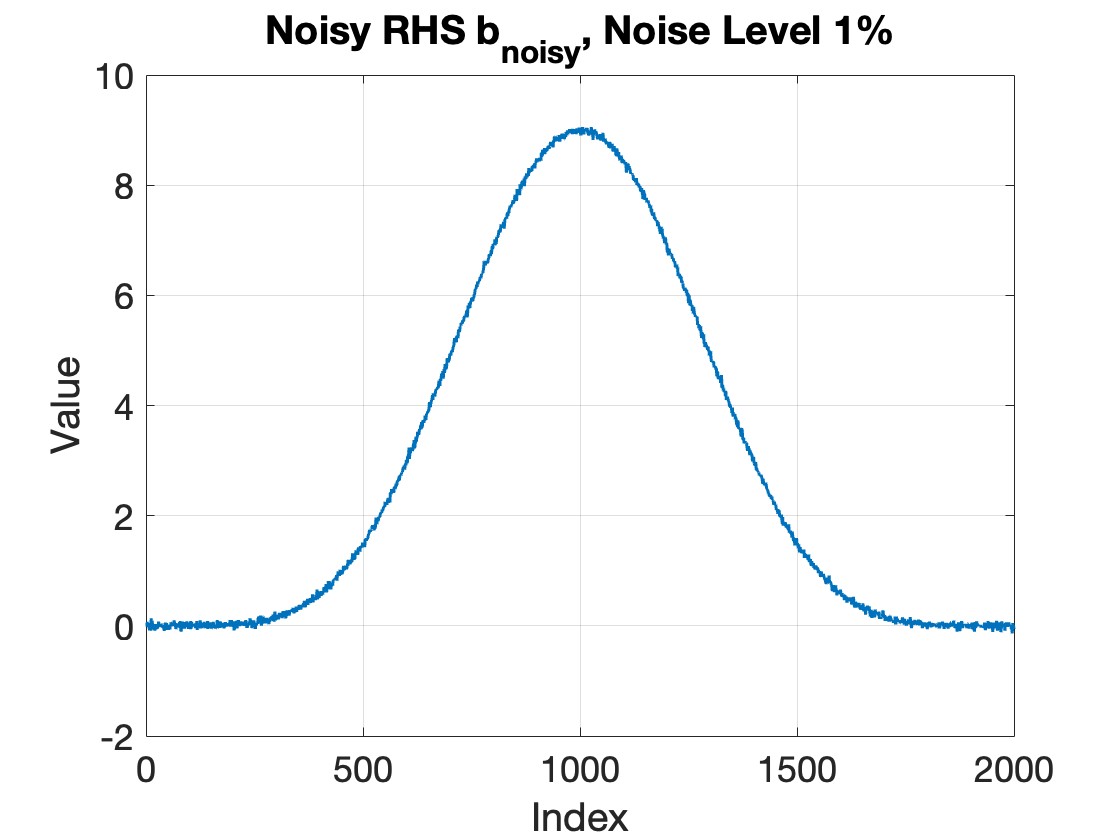}
    \caption{Visualization of the Phillips test problem. Left: True image \( x_{\text{true}} \). Middle: Right-hand side \( b \). Right: Noisy right-hand side \( b \).}
    \label{fig:test_problem}
\end{figure}

\begin{figure}[H]
    \centering
    \includegraphics[width=0.48\linewidth]{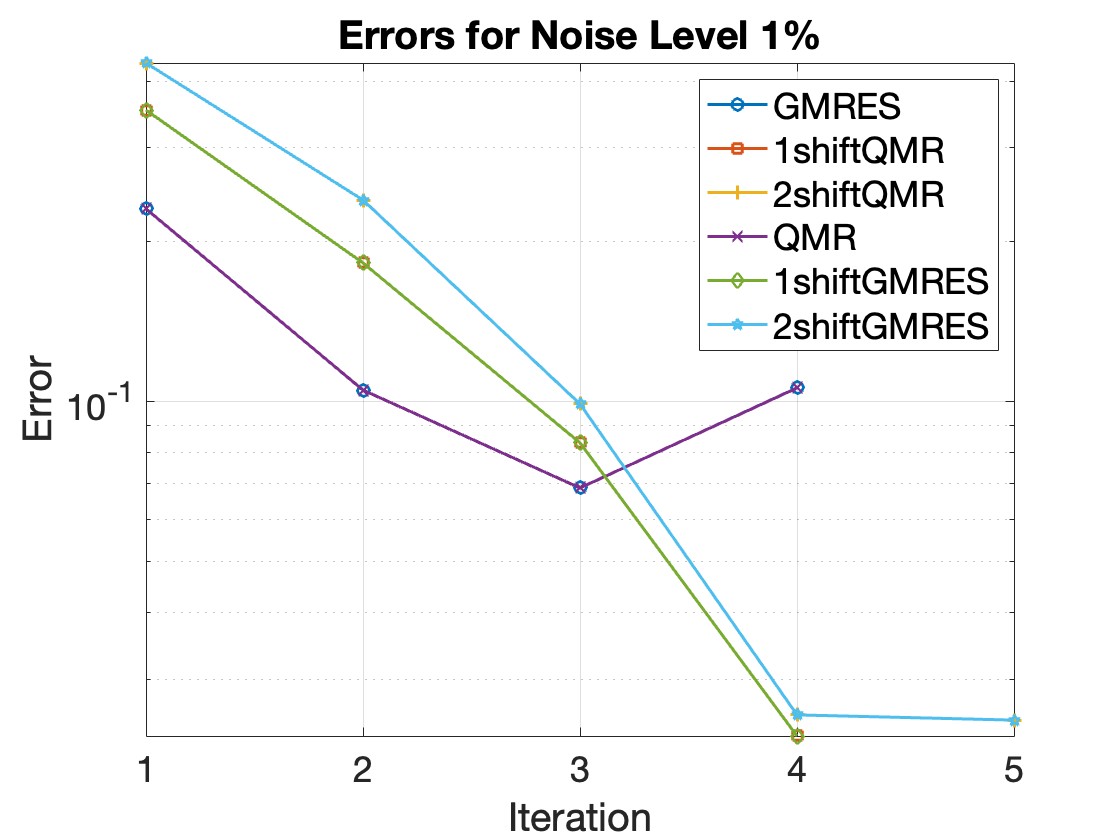}
    \includegraphics[width=0.48\linewidth]{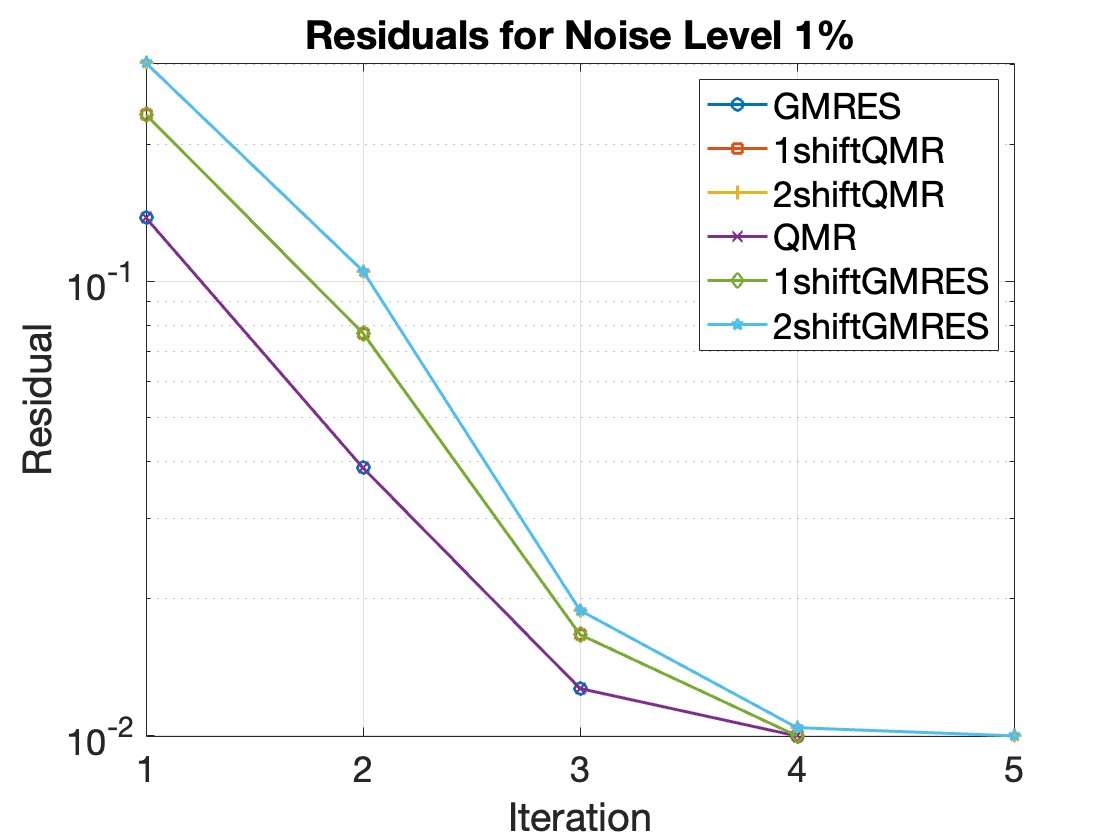}
    \caption{Comparison of relative reconstruction error for different methods. Left: Error comparison. Right: Residual comparison.}
    \label{fig:comparison}
\end{figure}

For this experiment, we set the size of the problem to be \( n = 2000 \) and used the Phillips problem with \( 1\% \) noise. The size of A is $2000 \times 2000$. From Figure 3, we observe that applying range restriction significantly improves the performance of both QMR and GMRES. While the residuals remain similar across methods, the error is noticeably lower when restriction is applied, implying a better recovery of the true solution. In the multi-shift case, 2-shift GMRES outperforms 1-shift GMRES, indicating that additional shifts can help further enhance solution accuracy. However, 2-shifted QMR does not show a notable improvement over 1-shift, suggesting that successive shifts may offer diminishing returns. The full table is given below.

\begin{table}[H]
    \caption{Relative reconstruction error for various solvers and noise levels.}
    \centering
    \scriptsize  
    \setlength{\tabcolsep}{3pt}  
    \renewcommand{\arraystretch}{0.95}  
    \begin{tabular}{|c|c|c|c|c|c|c|}
        \hline
        \textbf{Noise Level} & \textbf{GMRES} & \textbf{1-shift QMR} & \textbf{2-shift QMR} & \textbf{QMR} & \textbf{1-shift GMRES} & \textbf{2-shift GMRES} \\
        \hline
        0.1\% & 1.68e-02 & 9.91e-03 & 8.22e-03 & 1.64e-02 & 9.91e-03 & 8.22e-03 \\
        0.5\% & 5.79e-02 & 2.39e-02 & 2.50e-02 & 5.79e-02 & 2.39e-02 & 2.50e-02 \\
        1.0\% & 1.03e-01 & 2.51e-02 & 2.49e-02 & 1.03e-01 & 2.52e-02 & 2.49e-02 \\
        \hline
    \end{tabular}
    \label{tab:error_comparison}
\end{table}

\subsection{QMR shift comparison}

We analyze the effect of different shift values in QMR on the non-symmetric Philips problem. The graphs below correspond to the case with \( 5\% \) noise.

\begin{figure}[H]
    \centering
    \begin{minipage}{0.48\textwidth}
        \centering
        \includegraphics[width=\linewidth]{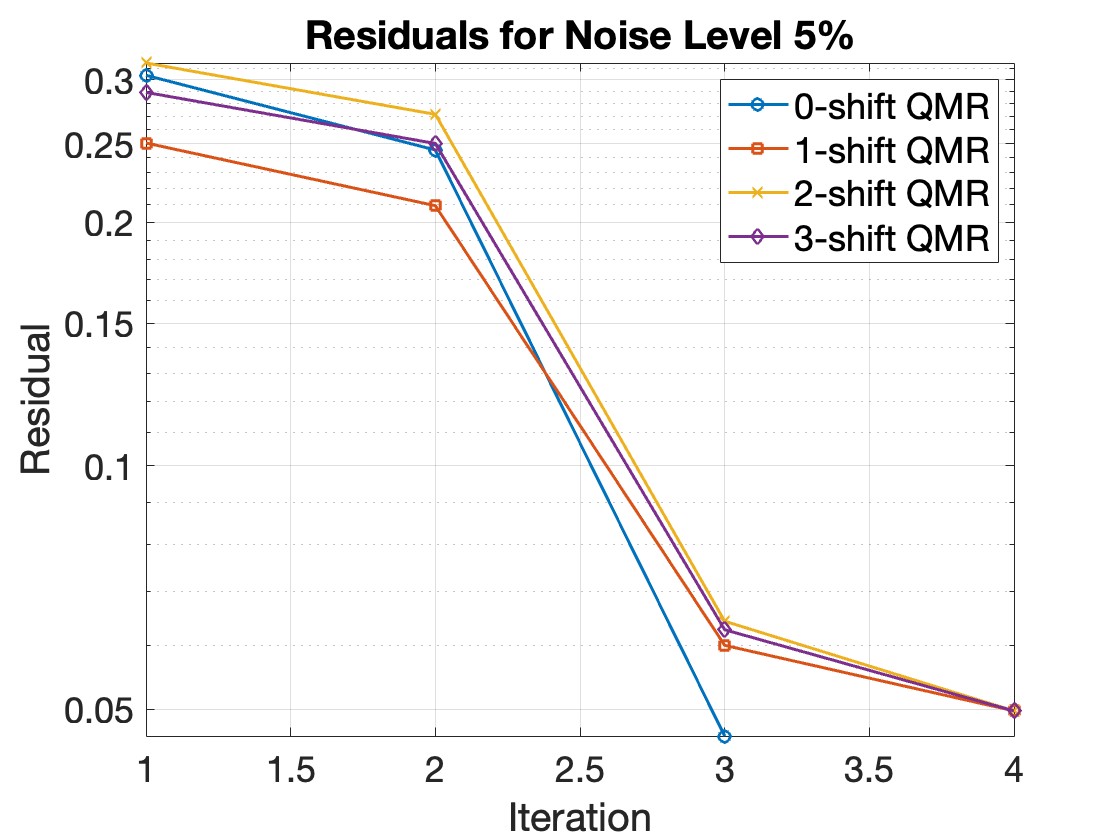}
        \label{fig:qmr_residual}
    \end{minipage}
    \hfill
    \begin{minipage}{0.48\textwidth}
        \centering
        \includegraphics[width=\linewidth]{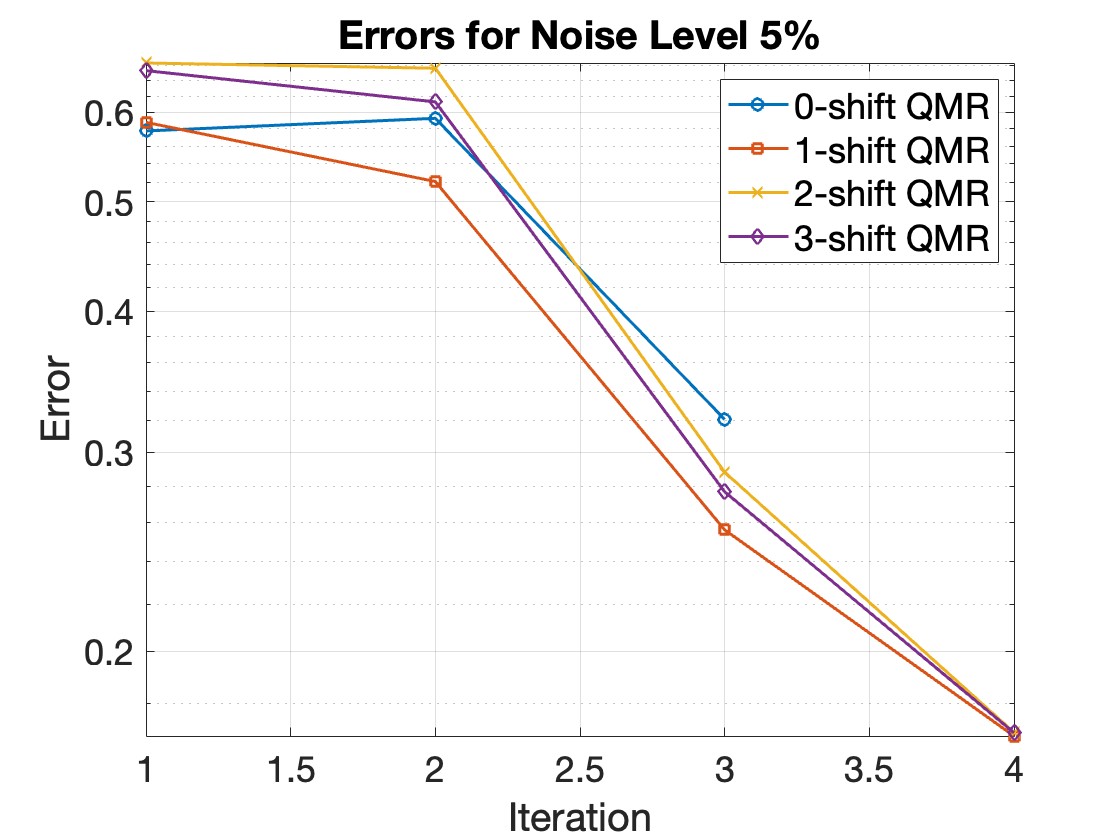}
        \label{fig:qmr_error}
    \end{minipage}
    \caption{Residual (left) and error (right) comparison of QMR with 0, 1, 2, and 3 shifts under 5\% noise.}
    \label{fig:qmr_shift_comparison}
\end{figure}

\begin{table}[H]
    \caption{Comparison of final errors for different QMR shifts.}
    \centering
    \footnotesize  
    \renewcommand{\arraystretch}{1.2}
    \begin{tabular}{|c|c|c|c|c|}
        \hline
        \textbf{Noise Level} & \textbf{standard QMR} & \textbf{1-shift QMR} & \textbf{2-shift QMR} & \textbf{3-shift QMR} \\
        \hline
        0.5\% & 1.748e-01 & 6.88e-02 & 5.10e-02 & 5.76e-02 \\
        1.0\% & 1.748e-01 & 1.15e-01 & 5.86e-02 & 4.97e-02 \\
        5.0\% & 3.209e-01 & 1.68e-01 & 1.70e-01 & 1.69e-01 \\
        \hline
    \end{tabular}
    \label{tab:qmr_shift_comparison}
\end{table}

We compare different shift values in QMR, ranging from 0 (standard QMR) to 1, 2, and 3 shifts. The final errors indicate that shifted QMR outperforms standard (0-shift) QMR. However, beyond the first shift, there is no significant improvement in error reduction.  

We observe that the 0-shift case has the lowest relative residual but also the highest relative error, demonstrating that a lower residual does not guarantee better reconstruction quality. Additionally, the number of shifts does not significantly impact QMR's performance, but increasing shifts requires more computation and storage. Without a specific reason to use multiple shifts, 1-shift QMR is likely the most practical choice for most situations.

\subsection{Performance under uncertain error norm bound}

 As discussed in the discrepancy principle section, a known noise level is required for it to function as optimally as possible. However, in practice, the true noise level may not always be accurately estimated. If the noise level is underestimated, the stopping criterion may not be triggered at the appropriate time, leading to over-iteration and potential amplification of noise. To illustrate this, we consider a scenario where the noise level is mistakenly assumed to be \( 0.01\% \), while the actual noise level is \( 1\% \). Without an accurate estimate, the discrepancy principle fails to stop the iteration properly, leading to excessive iterations and degraded solution quality.

\begin{figure}[H]
    \centering
    \begin{minipage}{0.48\textwidth}
        \centering
        \includegraphics[width=\linewidth]{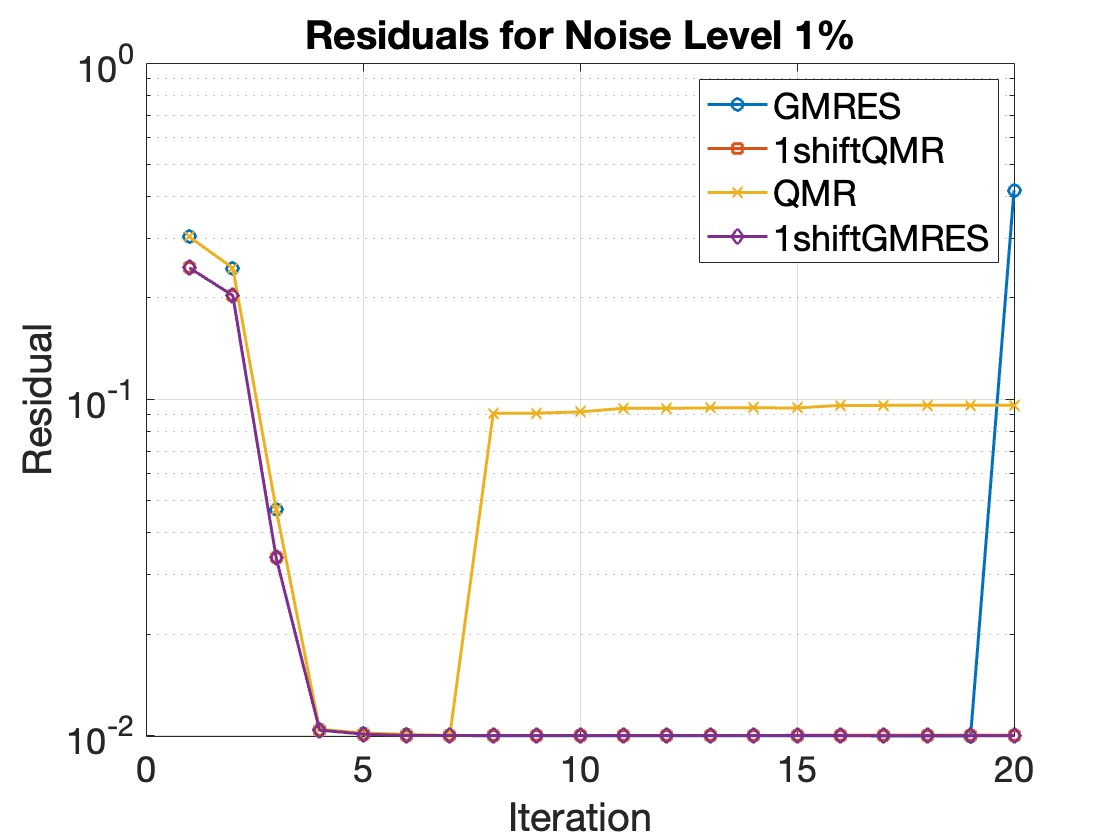}
        \label{fig:residual_no_discrepancy}
    \end{minipage}
    \hfill
    \begin{minipage}{0.48\textwidth}
        \centering
        \includegraphics[width=\linewidth]{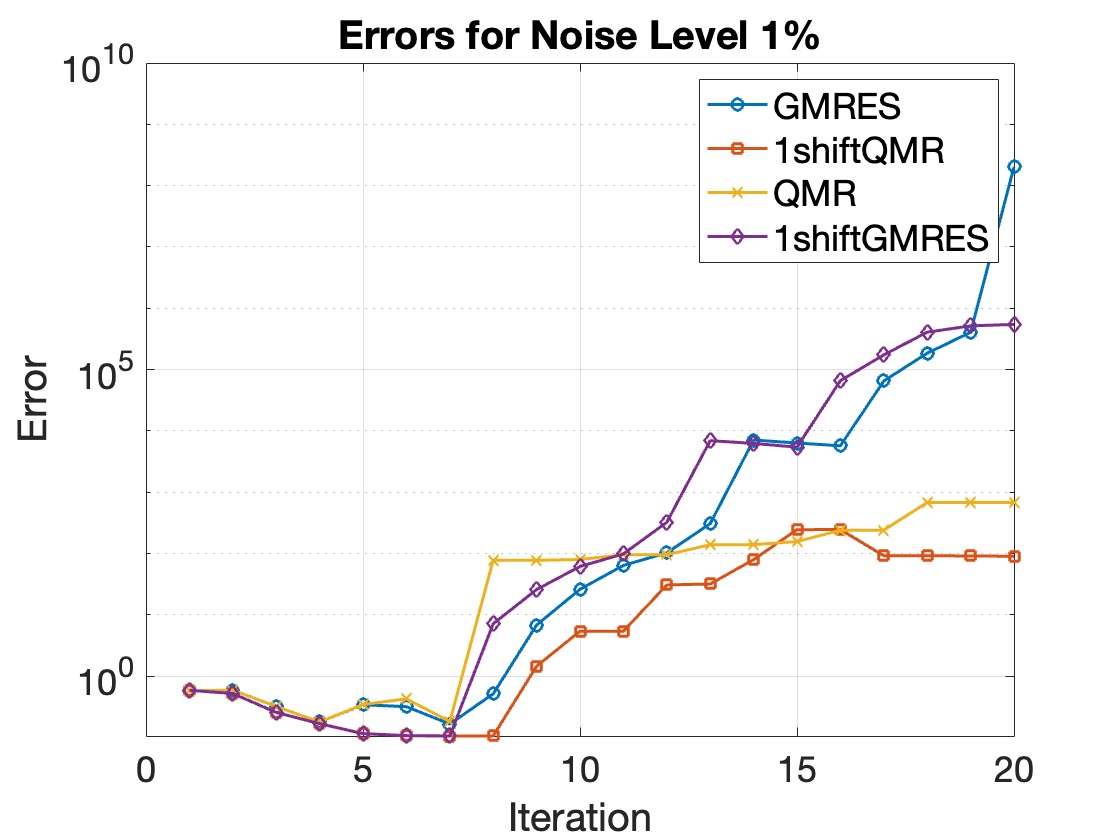}
        \label{fig:error_no_discrepancy}
    \end{minipage}
    \caption{Comparison of residual (left) and error (right) behavior for underestimate noise.}
    \label{fig:comparison_no_discrepancy}
\end{figure}

In the left of Figure 5, the error decreased first, and went up again. It increases as iterations progress, highlighting the semi-convergent nature of these solvers. The graph reveals that QMR and 1-shift QMR exhibit better semi-convergence behavior, achieving lower errors compared to GMRES and 1-shift GMRES, despite their residuals remaining similar.

Recall that in our methods, GMRES constructs an upper Hessenberg matrix \( H_m \) to approximate \( A \), while QMR builds a tridiagonal matrix \( T_m \) after $m$ steps. In the section of semiconvergence, we introduced the concept of the rapid decay of singular values amplifies errors in the solution. If the singular values of \( H_m \) or \( T_m \) decay more slowly, the impact of inverted noise may be less severe. To explore this idea, we compare the singular value decay of \( H_m \) and \( T_m \) for different shift.

\begin{figure}[H]
    \centering
    \includegraphics[width=0.6\linewidth]{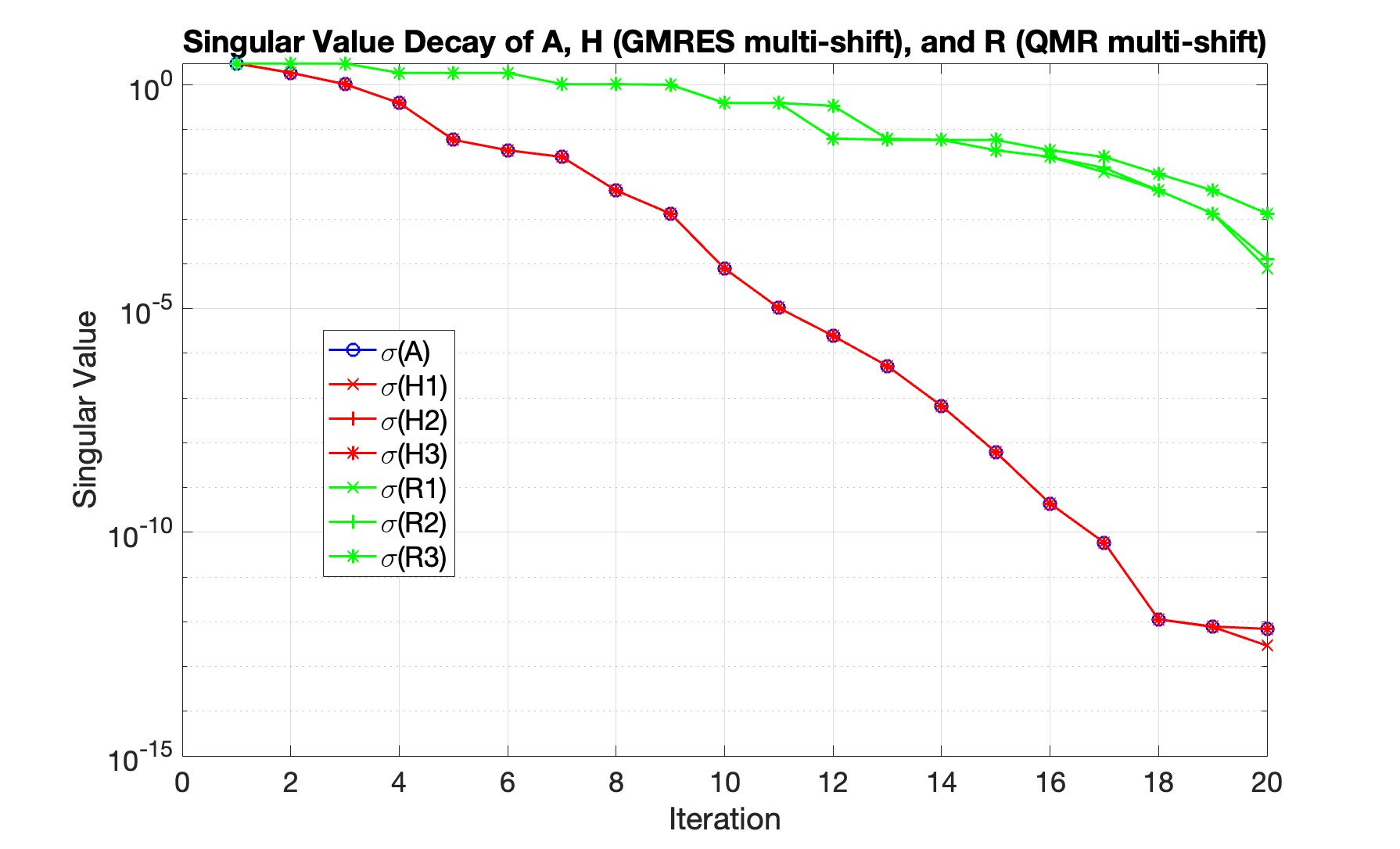}
    \caption{Comparison of Singular Value Decay in \( T_m \) from QMR and \( H_m \) from GMRES for different shift values.}
    \label{fig:singular_value_decay}
\end{figure}

From the graph above, we observe the singular values of QMR and GMRES for different shift values. Here, ``1'' represents no shift, ``2'' corresponds to one shift, and ``3'' denotes two shifts. For all GMRES and range restricted GMRES cases, we decompose the upper Hessenberg matrix \( H_m \), while for all QMR cases, we decompose the tridiagonal matrix \( T_m \). 

We note that in the range restricted QMR framework, we do not directly decompose \( T_m \), but rather use the upper triangular matrix \( R \) from the QR factorization of \( T_m \), \( T_m = Q R \), to approximate the spectral behavior of \( A \). We observe that the singular values of \( T_m \) decay significantly slower than those of \( H_m \) from GMRES. This aligns with our earlier discussion, where we established that both \( T_m \) and \( H_m \) serve as approximations to the original matrix \( A \). The slower decay of singular values in the QMR case suggests a reduced effect of inverted noise, which may explain the improved semi-convergence behavior of QMR.

In the situation of underestimating noise level, the range restricted QMR method performs better than both GMRES and the range restricted GMRES method. When the noise level is underestimated, iterative algorithms tend to over-iterate in an attempt to solve the problem. Under such conditions, the $\ell$-shifted QMR method may guarantee a lower error.

\subsection{2D problem: image deblurring}

In this experiment, we perform image deblurring with \( 1\% \) noise. The first image consists of six sub-images: the original image, the noisy image, and four reconstructions using different solvers.

\begin{figure}[H]
    \centering
    \includegraphics[width=0.8\linewidth]{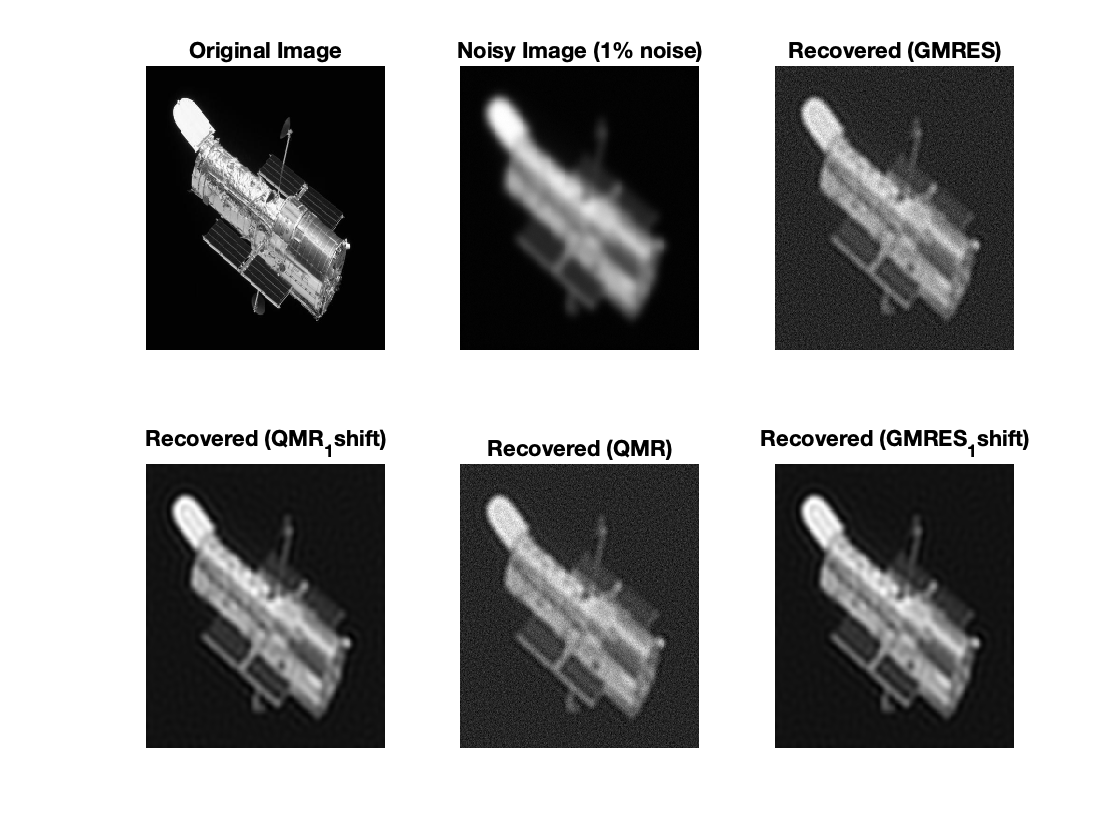}
    \caption{Comparison of image deblurring results. The first two images show the true and noisy images. The remaining four images represent reconstructions using GMRES, QMR, range restricted GMRES(1 shift), range restricted QMR(1 shift) .}
    \label{fig:image_deblurring}
\end{figure}

From this figure, we observe that the images recovered using 1-shift QMR and 1-shift GMRES are closer to the original image. These methods effectively balance sharpness and smoothness better than standard GMRES and QMR. To further analyze solver performance, we examine the residual and error plots.

\begin{figure}[H]
    \centering
    \begin{minipage}{0.48\textwidth}
        \centering
        \includegraphics[width=\linewidth]{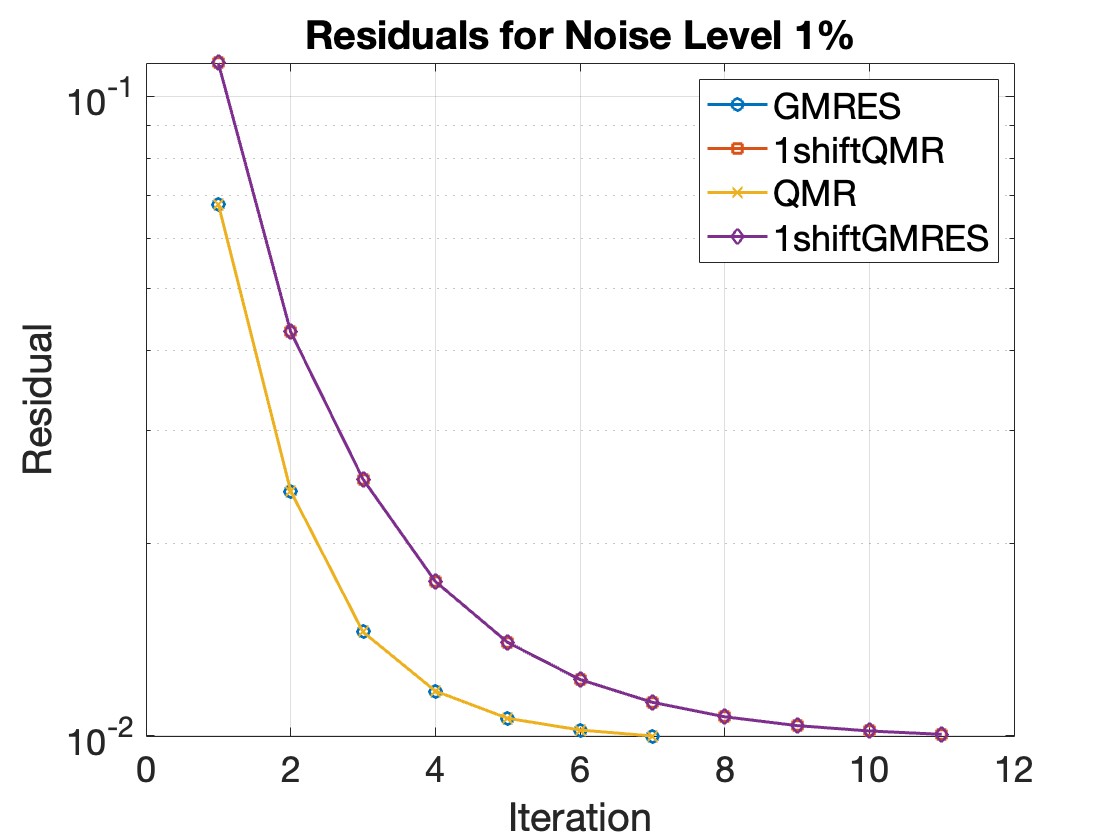}
        \label{fig:residual_deblurring}
    \end{minipage}
    \hfill
    \begin{minipage}{0.48\textwidth}
        \centering
        \includegraphics[width=\linewidth]{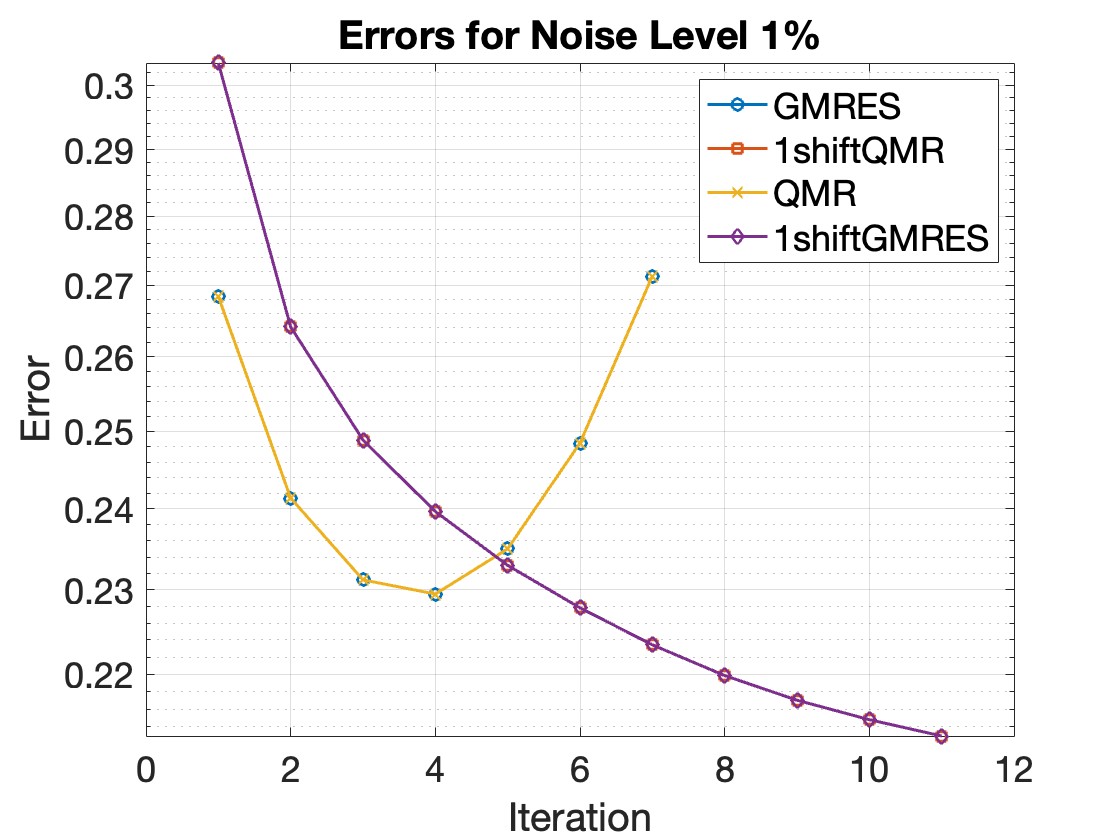}
        \label{fig:error_deblurring}
    \end{minipage}
    \caption{Comparison of residual (left) and error (right) behavior for different solvers in the image deblurring task.}
    \label{fig:deblurring_residual_error}
\end{figure}

The error plots confirm our visual observations. Range restricted QMR and GMRES consistently yield lower errors than their standard counterparts. However, when comparing shifted QMR and shifted GMRES, no significant difference in error reduction is observed. To quantify these observations, we summarize the final error values at different noise levels in the table below.

\begin{table}[H]
    \caption{Comparison of final errors for different solvers in image deblurring.}
    \centering
    \footnotesize  
    \renewcommand{\arraystretch}{1.2}
    \begin{tabular}{|c|c|c|c|c|}
        \hline
        \textbf{Noise Level} & \textbf{GMRES} & \textbf{1-shift QMR} & \textbf{QMR} & \textbf{1-shift GMRES} \\
        \hline
        0.5\% & 2.42e-01 & 2.05e-01 & 2.42e-01 & 2.05e-01 \\
        1.0\% & 2.71e-01 & 2.13e-01 & 2.71e-01 & 2.13e-01 \\
        5.0\% & 3.09e-01 & 2.34e-01 & 3.09e-01 & 2.34e-01 \\
        \hline
    \end{tabular}
    \label{tab:image_deblurring_errors}
\end{table}

The results indicate that 1-shift QMR and 1-shift GMRES produce better reconstructed images than their non-shifted counterparts. Shifted solvers consistently achieve lower errors compared to their unshifted versions, demonstrating the advantage of incorporating range restrictions in image deblurring. However, no significant difference is observed between shifted QMR and shifted GMRES in terms of final error, suggesting that both approaches benefit similarly from the shift technique. these experiments suggests that shifting improves solver performance, leading to better image reconstructions and reduced errors.

\section{Conclusions} \label{sec:conc}
In this work, we investigated Krylov subspace iterative methods for solving ill-posed inverse problems. Our primary focus was on the range restricted variants of GMRES and QMR methods. We then explored a range restricted variant of QMR and its comparison with range restricted GMRES. Our findings show that the range-restricted QMR method, which incorporates the range-restriction technique, outperforms both standard QMR and GMRES in solving ill-posed problems. A key advantage of range restricted QMR is its robustness in cases where noise levels are uncertain or underestimated. In such situations, non-restricted methods may continue iterating beyond the optimal stopping point, leading to the amplification of noise. The range restricted QMR addresses this issue by the nature of its tridiagonal matrix \(T_m\), thereby offering a more stable and reliable solution approach for ill-posed problems.


\bibliographystyle{siam}
\bibliography{references}

\end{document}